\documentclass[twoside,12pt, leqno]{amsart}
\usepackage{amsmath,amscd,amsthm,amssymb,amsxtra,latexsym,epsfig,epic,graphics}
\usepackage[matrix,arrow,curve]{xy}
\usepackage{graphicx}
\usepackage{diagrams}
\def\antiddot{\mathinner{\mkern1mu\raise1pt\vbox{\kern7pt\hbox{.}}\mkern2mu
        \raise4pt\hbox{.}\mkern2mu\raise7pt\hbox{.}\mkern1mu}}

\newcommand{\PP}{{\mathbb P}}
\newcommand{\QQ}{{\mathbb Q}}
\newcommand{\RR}{{\mathbb R}}

\newcommand{\HH}{{\rm{H}}}

\newcommand{\Ext}{{\rm{Ext}}}

\newcommand{\s}{\mathcal}

\newcommand{\cO}{{\s O}}

\newcommand{\punkt}{\hspace{-.3ex}\raise.15ex\hbox to1ex{\Huge.}}

\def \fix#1 {{\hfill\break \bf (( #1 ))\hfill\break}}

\DeclareMathOperator{\Sym}{Sym}

\DeclareMathOperator{\Hom}{Hom}

\DeclareMathOperator{\depth}{depth}

\DeclareMathOperator{\codim}{codim}

\newtheorem{theorem}{Theorem}[section]
\newtheorem{lemma}[theorem]{Lemma}
\newtheorem{proposition}[theorem]{Proposition}
\newtheorem{corollary}[theorem]{Corollary}
\newtheorem{conjecture}[theorem]{Conjecture}
\theoremstyle{definition}

\newtheorem{question}[theorem]{Question}

\newtheorem{remark}[theorem]{Remark}

\newtheorem{example}[theorem]{Example}

\def\PP{{\mathbb P}}
\def\P{{\mathbb P}}
\def\QQ{{\mathbb Q}}

\def\fix#1{{\bf ***Fix:} #1 {\bf ***}}

\def\Rbar{{\overline R}}
\def\Ibar{{\overline I}}
\def\mm{{\frak m}}
\def\RR{{\mathcal R}}

\def\tR{{\tilde R}}
\def\tI{{\tilde I}}
\def\tJ{{\tilde J}}
\def\tK{{\tilde K}}

\newarrow{Iso} -----

\def\Rbar{{\overline R}}
\def\Ibar{{\overline I}}

\def\abar{{\overline \alpha}}

\def\abar{{\overline a}}
\def\End{{\rm End}}

\def\red{{\rm red}}

\def\lbracket{{[\kern-1.5pt[}}
\def\rbracket{{]\kern-1.5pt]}}

\newtheorem{proposition-definition}[theorem]{Proposition-Definition}

\makeatletter
\def\Ddots{\mathinner{\mkern1mu\raise\p@
\vbox{\kern7\p@\hbox{.}}\mkern2mu
\raise4\p@\hbox{.}\mkern2mu\raise7\p@\hbox{.}\mkern1mu}}
\makeatother

\usepackage{times}
\newdimen\x \x=12pt

\usepackage{color}

\usepackage{color}
\usepackage[usenames,dvipsnames,svgnames,table]{xcolor}

\usepackage[breaklinks,bookmarksopen,bookmarksnumbered,urlcolor=blue]{hyperref}
\hypersetup{colorlinks=true,backref=true,citecolor=blue}

\author{David Eisenbud, Craig Huneke, and Bernd Ulrich
}

\title{Residual Intersections and Linear Powers}
\begin{document}

\begin{abstract} If $I$ is an ideal in a Gorenstein ring $S$, and $S/I$ is Cohen-Macaulay, then the same is true for any 
linked ideal $I'$; but such statements hold for residual intersections of higher codimension only
under restrictive hypotheses, not satisfied even by ideals as simple as the ideal $L_{n}$ of minors of a generic $2 \times n$ matrix when $n>3$. 

In this paper we initiate the study of a different sort of Cohen-Macaulay property that holds for certain general residual intersections of the maximal (interesting) codimension, one less than the analytic spread of $I$. For example, suppose that $K$ is the residual intersection of $L_{n}$ by $2n-4$ general quadratic forms in $L_{n}$. In this situation we analyze $S/K$ and show that 
$I^{n-3}(S/K)$ is a self-dual maximal Cohen-Macaulay $S/K$-module with linear free resolution over $S$. 

The technical heart of the paper is a result about ideals of analytic spread 1 whose high powers are linearly presented.
\end{abstract}

\maketitle

\section*{Introduction}
\let\thefootnote\relax\footnote{
\noindent AMS Subject Classification:\\
Primary: 13C40, 13H10, 14M06, 14M10;
Secondary: 13D02 , 13N05, 14B12, 14M12.\smallbreak
We are grateful to the NSF for partial support, and to MSRI and the Polini household for hosting us during some of the work on this paper.}
\bibliographystyle{alpha}

Let $S = k[x_1,\dots,x_d]$ be a standard graded polynomial ring over a field, and let $I\subset S$ be
an ideal of codimension $g$. Let $J = (a_1,\dots, a_s)$ be an ideal generated by $s$ 
forms in $I$, and consider the residual ideal $K := J:I$. If $\codim K\geq s$ then  $K$ (or $S/K$) is 
said to be the \emph{$s$-residual intersection of $I$ with respect to $J$}, and it is called \emph{geometric}
if, in addition, $\codim (I+K) >s$.

Under  strong hypotheses the residual intersection $R = S/K$ is Cohen-Macaulay and, up to shifts, the canonical module $\omega_R$
of $R$ is isomorphic to $I^{s-g+1}R$  and  $I^{j}R$ is $\omega_R$-dual to $I^{s-g+1-j}R$ for $0\leq j\leq s-g+1$ 
(see~\cite{EU} for a summary of the situation). For example, these conclusions are true
when $I$ is the ideal of maximal minors of a sufficiently general $(n-1)\times n$ matrix over a polynomial ring.

However, all these things fail even when $I$ is generated by the maximal minors of a generic $2\times n$ matrix with $n\geq 4$. The main contribution of this paper is to construct a natural rank 1, self-dual, maximal Cohen-Macaulay module over 
 certain residual intersections of such ideals and many others. The following special case of our main result, Theorem~\ref{main}, will convey the flavor:

\begin{theorem}\label{special case of 2.3}
Let $S$ be a standard graded polynomial ring over a field of characteristic 0, let $I\subset S$ be
a homogeneous ideal generated in a single degree $\delta$ with analytic spread $\ell$, and let  $J\subset I$ be generated by $\ell-1$ general elements of degree $\delta$. Suppose that $J:I$ is an $(\ell-1)$-residual intersection of $I$, and set $R = S/(J:I)$. If all sufficiently high powers of $I$ are
linearly presented then, for all $\rho\gg 0$,
$
I^{\rho}R
$
 is a maximal Cohen-Macaulay $R$-module with linear resolution as an $S$-module, and is, up to shift, $\omega_{R}$-self dual.
\end{theorem}

The restriction to $(\ell-1)$-residual intersections is natural because in that case $IR$ has analytic spread 1, so all high powers 
are isomorphic. Thus it makes sense to speak of their asymptotic structure.

The idea of the proof is to reduce to the case of analytic spread 1, and then use the fact that $I^\rho R$, for large $\rho$, can be given the structure of a commutative algebra. The key point is to show that the condition of linear presentation is preserved in the reduction.

Theorem~\ref{special case of 2.3}  applies, in particular, to the ideal of maximal minors of a generic $2\times n$
matrix, and more can be said in that case:

\begin{theorem}\label{2xn0} Let $S= k[x_{1},\dots,x_{n}, y_{1},\dots,y_{n}]$ be a polynomial ring over
  a field $k$ of characteristic 0. Suppose that $I \subset S$  is the ideal of $2\times 2$ minors of the generic matrix
  $$
\begin{pmatrix}
 x_{1}&\dots&x_{n}\\y_{1}&\dots&y_{n}
\end{pmatrix}.
$$
 The ideal $I$ has analytic spread $\ell = 2n-3$.
 
Let $J$ be an ideal
  generated by $\ell-1$ general quadrics in $I$. The ring $R = S/(J:I)$  is unmixed of codimension $2n-4$, with isolated singularity.
  If $k$ is algebraically closed, then $J:I$ is the intersection of $\frac{1}{n-1}{2(n-2)\choose n-2}$
linear prime ideals. 

 The canonical module $\omega_{R}$ is isomorphic to $I^{n-3}R(2n-10)$. Furthermore, for all  $\rho \geq n-3$, the module $I^\rho R$  
  is isomorphic up to shift to $\omega_R$, and is a maximal Cohen-Macaulay $R$-module with linear resolution
  as an $S$-module.  \end{theorem}

Examples suggest that part of the conclusion of Theorem~\ref{special case of 2.3} holds much more generally:

 \begin{conjecture}\label{main conjecture0} 
 Let $S$ be a standard graded polynomial ring over an infinite field, let $I\subset S$ be
a homogeneous ideal generated in a single degree $\delta$ with analytic spread $\ell$, and let  $J\subset I$ be generated by $\ell-1$ general elements of degree $\delta$. Set $R = S/(J:I)$. If $I$ is unmixed and $J:I$ is a geometric $(\ell-1)$-residual intersection of $I$ then
$
I^{\rho}R
$
 is a maximal Cohen-Macaulay $R$-module for all $\rho\gg 0$. 
 \end{conjecture}

Further examples suggest that Conjecture~\ref{main conjecture0}  holds independently of the assumption that $R$ is obtained as a residual intersection:

\begin{conjecture}\label{analytic spread 1}
Let $R$ be a standard graded algebra over a field. Assume that $R$ is reduced and equidimensional,   and that $\omega_{R}$ is generated in 
a single degree. If $\omega_{R}$ has analytic spread 1, then some power of $\omega_{R}$ (as a fractional ideal) is $\omega_{R}$-self dual up to a shift.
\end{conjecture}

Examples we have seen also support another version:

\begin{question}\label{analytic spread 1 with truncation}
Let $R$ be a standard graded  ring over a field, and let $I\subset R$ be a homogeneous ideal of positive codimension. 
Suppose that:
\begin{enumerate}
 \item $R$ is reduced and the truncation $I'$ of $I$ in the degree of the highest generator of $I$ has analytic spread 1;
 \item Some power of $I$ is isomorphic to a power of a canonical ideal of $R$ up to a shift.
\end{enumerate}
Are the high powers of $I'$ always $\omega_R$-self-dual modules? 
\end{question}

We shall see in Section~\ref{S examples} that none of the  hypotheses can be dropped.
Conjecture~\ref{main conjecture0} and has surprising consequences:

\begin{proposition}\label{cor of conjecture}
 Suppose that $S = k[x_1,\dots,x_d]$ is a standard graded polynomial ring, and $I\subset S$ is an ideal generated by forms of a single degree. If Conjecture~\ref{main conjecture0} is true then: 
\begin{enumerate}
 \item Suppose that $I$ is generically a complete intersection of codimension $g$ generated by $g+1$ elements. If $I$ is unmixed then $S/I$ is Cohen-Macaulay.
 \item If $I\subset S$ is an unmixed ideal of linear type (that is, the Rees algebra $\mathcal R(I)$ of $I$ is equal to the symmetric algebra of $I$) then 
 $$
 \mathcal R(I)_{(x_1,\dots,x_d) \mathcal R(I)}
 $$
   is Cohen-Macaulay.
 \end{enumerate}
\end{proposition}
The assertion of item (1) is known to be true, independent of the conjecture, both when $I$ has codimension 2 \cite[Theorem 2.1]{EG}, and also when $S/I$ is equidimensional and locally of depth $\geq \frac{1}{2}\dim (S/I)_P$ at every prime $P$ containing $I$ (Theorem~\ref{ACI} below). In Section~\ref{proof of 0.5} we  prove Proposition~\ref{cor of conjecture}, and note two special cases where the conjectures are verified. 

\vspace{.1in}

 \section{The module of interest}\label{S basic properties}
 
 The analytic spread of an ideal plays an important role in this theory. 
If $R$ is a positively graded algebra over a field $k$, with maximal homogeneous ideal $\mm$, and
$I\subset R$ is an ideal generated by forms of a single degree $\delta$,
then the analytic spread may be defined as:
  $$
 \ell(I) := \dim k[I_\delta] = \dim (k \oplus I/\mm I \oplus I^2/\mm I^2\oplus \cdots).
 $$
Assuming that $k$ is infinite, $\ell(I)$ is the smallest number of generators of a homogeneous ideal over which $I$ is integral, and such an ideal may be taken to be the ideal generated by $\ell(I)$ general forms in $I$ of  degree $\delta$. The reduction 
number $r(I)$ of $I$ is the smallest integer $r \geq 0$ so that $I^{r+1}=JI^r$ for some homogeneous $\ell(I)$-generated ideal
$J$ over which $I$ is integral.

If $s<\ell(I)$  and $K$ is an $s$-residual intersection of $I$ with respect to $J$, so that $\codim K \geq s$, then $I$ cannot be integral over $J$, 
so if $R$ is equidimensional then, by \cite[4.1]{McAdam}, $\codim K = s$. However, this may not be the case when $s\geq \ell(I)$. 

Now assume that $k$ is infinite. When $s<\ell(I)$ and $J$ is generated by $s$ general forms of degree $\delta$ in $I$, Proposition~\ref{basic setup higher ell} implies that
the modules
$I^\rho (S/K)$ are nonzero. These, for $\rho\gg 0$ and $s = \ell-1$, are the modules that are of interest to us.

\begin{proposition-definition}\label{basic setup ell=1}
Suppose that $R$ is a standard graded algebra over an infinite field $k$,  that $I\subset R$ is generated by forms of a single degree $\delta$, and that $I$ has analytic spread
$\leq 1$. Let $\Rbar := R/(0:_RI^\infty)$ and $\Ibar := I\Rbar$.
\begin{enumerate}
 \item For $\rho \gg 0$ the module $ I^\rho (\rho \delta)$ is, up to homogeneous isomorphism, independent of $\rho$ and
$I^\rho$ maps isomorphically to $\Ibar^\rho$.
 If $I$ contains a nonzerodivisor of $R$, then this is true for all $\rho$ greater than or equal to the reduction number $r(I)$.

\end{enumerate}
\smallbreak
Let $M := M(I)$  be the stable value of $ I^\rho(\rho \delta)$. 
\smallbreak
\begin{enumerate}\setcounter{enumi}{1}
\item  $E := \End_R(M)\cong M$ as  graded $R$-modules. 

\item Let $a\in I$ be a general homogeneous element of degree $\delta$. Write
$\Ibar$ and $\abar$ for the images of $I$ and $a$ in $\Rbar$. The element $\abar$ is a nonzerodivisor on $\Rbar$,  and 
$$
E \cong \Rbar[(1/\abar)\Ibar] = (1/\abar^\rho)\Ibar^\rho ,
$$
 which is
the coordinate ring of the blowup of $\Rbar$ along $\Ibar$.

\item If $R$ is reduced away from $V(I)$, then
 $$
 \Rbar = R_\red/(0:(IR_\red)^\infty) = R/(\sqrt 0:I^\infty)
 $$
is reduced.
\end{enumerate}

\end{proposition-definition}

\begin{proof}
\noindent{\bf  (1)} Let $\epsilon$ be an integer so large that $0:I^\infty = 0:I^\epsilon$. 
By the Artin-Rees Lemma,
$I^\rho  \cap (0:_R I^\epsilon ) \subset I^\epsilon \, (0:_R I^\epsilon ) =0$ for large $\rho$. Thus  $I^\rho$ maps
isomorphically to $\Ibar^\rho$. Note that
 $\Ibar$ has positive grade; in particular
 $\abar$ is a nonzerodivisor on $\Rbar$. 

Because $I$ has analytic spread $\leq 1$, an element $\abar$ as in $(3)$  generates a reduction of $\Ibar$. Thus the
multiplication map
$$
\abar: \Ibar^\rho(\rho \delta) \to \Ibar^{\rho+1}((\rho+1) \delta)
$$
 is surjective, and since $\abar$ is a nonzerodivisor on $\Rbar$ the map is injective as well.
 This proves that  $M$ is independent of $\rho\gg 0$. 
\smallbreak
\noindent{\bf  (2)}  We have $\End_R(M) = \End_\Rbar(M)$. Writing $Q$ for the total ring of
quotients of $\Rbar$,  we have
$$
\Ibar^\rho = \abar^\rho \Ibar^\rho:_Q \abar^\rho \Rbar \supset \abar^\rho \Ibar^\rho:_Q \Ibar^\rho \supset \Ibar^\rho
$$
because  $\abar \Rbar$ is a reduction of $\Ibar$. Thus $\Ibar^\rho =\abar^\rho \Ibar^\rho:_Q \Ibar^\rho$. Further $\abar^\rho \Ibar^\rho:_Q \Ibar^\rho \cong \End_\Rbar(\Ibar^\rho)(-\rho\delta)$ proving the assertion.

\smallbreak
\noindent{\bf  (3)} Because $\abar\Rbar$ is a principal reduction of $\Ibar$, the blowup has only one affine chart,
$\Rbar[(1/\abar)\Ibar ]$. Further $(1/\abar)\Ibar \subset (1/\abar^{2})\Ibar^2\subset \cdots$, and the sequence
of fractional ideals stabilizes. Thus for $\rho\gg 0$ we have
$$
\Rbar[(1/\abar)\Ibar] = (1/\abar^{\rho})\Ibar^\rho = \Ibar^\rho:_Q \Ibar^\rho \cong E.
$$

\smallbreak
\noindent{\bf  (4)} If $R$ is reduced away from $V(I)$ then $I^\gamma \sqrt {0} = 0$ for some $\gamma$. By the Artin-Rees Lemma,
$I^\rho \cap \sqrt{0}\subset I^\gamma\sqrt {0} = 0$ for $\rho \gg 0$. Thus, for $\rho \gg 0$,
$$
\sqrt 0:I^\rho = (I^\rho\cap \sqrt 0):I^\rho = 0:I^\rho.
$$
 \end{proof}

 The next result provides a different description of $I^\rho \Rbar$ in a general setting:
 
\begin{proposition}\label{basic setup higher ell}
Let $S$ be a Noetherian algebra over an infinite field  $k$, and let $I$ be an ideal of $S$. 
Let $J$ be an ideal generated by a sequence of general $k$-linear combinations of generators of $I$. Let $R = S/(J:I)$, and let $\Rbar = S/(J:I^\infty)$.

The natural surjection
$$
I^\rho/JI^{\rho-1} \to I^\rho \Rbar
$$
is an isomorphism for  $\rho\gg 0$. 
\end{proposition}
 
 \begin{proof} Let $\epsilon$ be an integer so large that $H: = J:I^\infty = J:I^\epsilon$. 
For $\rho \gg 0$  the Artin-Rees Lemma gives
$I^\rho \cap H \subset I^\epsilon H \subset J$. 
Furthermore, the generators of $J$ are a superficial sequence for $I$, so \cite[Lemma 8.5.11]{SH}
gives
$
I^\rho \cap J =JI^{\rho-1}.
$
Thus
$I^\rho \cap H=JI^{\rho-1}$
as required.
\end{proof}

Suppose that $S = k[x_1,\dots, x_d]$ is a standard graded polynomial ring over an infinite field $k$, and let 
$I$ be an ideal generated by forms of a single degree $\delta$. If the analytic spread of $I$ is $\ell$ and $J$ is
generated by $\ell-1$ general forms of degree $\delta$ then, setting $R = S/(J:I)$, the ideal $IR$ has analytic spread $\leq 1$,
so  we may apply Proposition-Definition~\ref{basic setup ell=1}. In this case the module $M= M(IR)$ can be expressed as $(I^\rho/JI^{\rho-1})(\rho\delta)$, as we see from Proposition~\ref{basic setup higher ell}.

\vspace{.05in}
 \section{Ideals whose powers have linear presentation}
 
 The following is our main general result.
 
 \begin{theorem}\label{main}
 Let $S$ be a standard graded polynomial ring over an infinite field $k$. Let $0\neq I\subset S$ be an ideal generated by forms of
 a single degree $\delta$. Let
 $\ell = \ell(I)$ be the analytic spread of $I$. Let $J\subset I$ be generated by $\ell-1$ general homogeneous elements of degree $\delta$.
Set 
$$
R= S/(J:I), \ \Rbar = S/(J:I^\infty), \ \Ibar = I\Rbar.
$$ 
Let $\abar\in \Ibar$ be a general form of degree $\delta$ and let
$M = M(IR)$ be as in 
Proposition-Definition~\ref{basic setup ell=1}.

If $R$ is reduced away from $V(I)$ and all sufficiently high powers of $I$ are linearly presented,
then 
\begin{enumerate}
  
 \item\label{2.1.1} $\Rbar$ is equidimensional,  $M$ is a maximal Cohen-Macaulay $\Rbar$-module and an Ulrich module, and
 $M$ is $\omega_{\Rbar}$-self-dual.

 \item\label{2.1.2}  $\dim \Rbar = \dim R$ if and only if  $J:I$ is an $(\ell - 1)$-residual intersection of $I$, that is, $\codim (J:I) \geq \ell - 1$.
 In this case,  $M$ is a maximal Cohen-Macaulay $R$-module and an Ulrich module, and
 $M$ is $\omega_{R}$-self-dual.
 
  \item \label{2.1.3} As a graded $R$-module, $\End_R(M) \cong M = I^\rho R(\rho\delta) \cong \Ibar^{\rho}(\rho\delta)$ for $\rho\gg 0$.
  As a graded ring, $\End_R(M)$ is isomorphic to the blow-up 
$$
\tilde R := \Rbar[\abar^{-1}\Ibar]
$$
 of $\Rbar$ along $\Ibar$ and is regular. 
 In particular $\tilde R$ is the normalization of $\Rbar$, and the conductor of $\Rbar$  is $\Rbar \abar^{\rho}:_{\Rbar}\Ibar^{\rho}$.

\end{enumerate}\end{theorem}

\begin{remark}\label{Bertini}
If the characteristic of $k$ is 0, then  $R$ is automatically reduced away from $V(I)$ by Bertini's Theorem
 (\cite{Flenner}).\end{remark}

\begin{remark} If $I$ satisfies $G_{\ell-1}$, then the ideal $K$ in Theorem~\ref{main} is an $(\ell-1)$-residual intersection.
 
\end{remark}
 
\begin{theorem}[Examples]\label{list}
 The following classes of ideals in a polynomial ring in $d$ variables over a field of characteristic 0 all satisfy the hypotheses and conclusions of parts (1)-(3) of Theorem~\ref{main}:
 
\begin{enumerate}
\item ideals of $m\times m$ minors of  $m\times n$ matrices $A$ of linear forms such that
either $\codim I_{m}(A) =d$ or
$$
\codim I_{k}(A) \geq \min((m-k+1)(n-m)+1, d)\hbox{ for } 2\leq k\leq m;
$$
\item strongly stable monomial ideals generated in one degree;
\item products of ideals of linear forms;
\item polymatroidal ideals;
\item monomial ideals generated in degree 2 and having linear resolution;
\item linearly presented ideals of dimension 0, and ideals of dimension 1 that have linear resolutions for the first
$\lceil (d -1)/2\rceil$ steps.
\item  truncations $I = I'_{\geq t}$ of  homogeneous ideals $I'$ at degree $t$ if $I'$ is generated in degrees $\leq t$ and
$I'$ has a homogeneous reduction generated in degrees $\leq t-1$.
\item linearly presented ideals of fiber type, such as linearly presented ideals satisfying $G_{d}$ that are perfect  of codimension 2 or Gorenstein  of codimension 3.
\end{enumerate}
\end{theorem}

More precisely, in the first 5 cases every power of the ideal in question actually has a linear resolution; in cases 6 and 7 all large powers
have linear resolution; and in case 8 all powers are linearly presented.\\

\begin{proof} 
(1): See Theorem~\ref{Bruns et al}, the second case of which is the main theorem of \cite{BCV}.\\
(2)--(4): These assertions are copied from the list in \cite{BCV}, and were proven in \cite{CH}. \\
(5): See \cite[Theorem 3.2]{HHZ}.\\
(6): See \cite[Theorem 7.1 and Corollary 7.7]{EHU}. \\
(7): See Proposition~\ref{trunc}.\\
(8): Symmetric powers of linearly presented ideals are always linearly presented; and fiber type implies that the additional relations
on the generators of the powers all have degree 0.  The given classes of ideals are of this type by
\cite{MU} and \cite{KPU}, respectively. In the case of perfect codimension 2 ideals, all the powers have linear resolution by item (1) and the Hilbert-Burch Theorem.
\end{proof}
 
\begin{proposition}\label{trunc}
Suppose that $I'$ is a homogeneous ideal in a standard graded polynomial ring. Suppose further that $I'$ is generated in degrees $\leq \delta$ and that $I'$ has a homogeneous reduction generated in degrees $\leq a$. 

If $t\geq \max(\delta, a+1)$, then the high powers of the truncated ideal $I'_{\geq t}$
all have linear resolution.
\end{proposition}

\begin{proof} 
If $t\geq \delta$  then $(I'_{\geq t})^\rho = (I'^\rho)_{\geq t\rho}$. 
Moreover, for large $\rho$ the regularity of
$I'^\rho$ grows as a linear function bounded by $a\rho+b$, see \cite{K} or \cite{CHT}.

Thus if   $t\geq \max(\delta, a+1)$ then $(I'_{\geq t})^\rho = (I'^{\rho})_{\geq t\rho}$ since $t\geq \delta$, and 
$ I'^{\rho}$ has regularity $\leq t\rho$ since $t\geq a+1$.
\end{proof}

\begin{remark}\label{linearly presented examples}
We do not know of an ideal whose high powers have linear presentation but not linear resolution.
Linearly presented ideals of fiber type may provide such examples.
\end{remark}

The following consequence of Theorem~\ref{main} gives some evidence for a positive answer to Question~\ref{analytic spread 1 with truncation}.

\begin{corollary}
Suppose that $I$ is a homogeneous ideal in a polynomial ring $S$ in $d$ variables over an infinite perfect field, and let
 $t$ be such that $I$ is generated in degree $\leq t$ and the ideal $I'^\rho$ has linear presentation for all large $\rho$  (in char p we assume $t >$ the generator degree of $I$). Let $J$ be an ideal generated by
 $d-1$ general forms of degree $t$ in $I'$.
 Set $K = J:I$ and $R := S/K$. 

 The conclusions of Theorem~\ref{main} hold
 if we replace $\delta$ by $t$, $IR$ by $I'R$, $\Ibar = I\Rbar$ by $I'\Rbar$;  and  $M(IR)$ by $M(I'R)$. 
 \end{corollary}

\begin{proof}
 We may assume that $\ell(I') = d$, since otherwise $I'^\rho \subset J$ for $\rho\gg 0$.
 Note that $J:I'^\infty = J:I^\infty$, so $\Rbar = S/(J:I'^\infty)$.  By~\cite{Flenner}, $R$ is reduced away from $V(I') = V(I)$.
 We now apply Theorem~\ref{main} with $I := I'$.  
\end{proof}

To prove Theorem~\ref{main} we reduce to the case $\ell=1$ by factoring out a general $(\ell-2)$-residual intersection and proving that the hypothesis of linearly presented powers is preserved. We then use the following more general result:

\begin{theorem}\label{stable residuals and linearly presented powers}
Let $S$ be a standard graded polynomial ring over an infinite field,  let $R$ be a homogeneous factor ring of $S$, and let $I\subset R$ be an ideal. 

Suppose that:
\begin{itemize}
 \item $I$ is generated by forms of a single degree $\delta$ and has analytic spread $\leq 1$. 
 \item $R$ is reduced away from $V(I)$.
 \item All sufficiently large powers of $I$ are linearly presented as $S$-modules.
\end{itemize}

Set
$\Rbar = R/(0:I^\infty)$ and  $\Ibar = I\Rbar$, and let $\abar\in \Ibar$ be a general element of degree $\delta$,
and let 
$$
\tilde R := \Rbar[\abar^{-1}\Ibar]
$$
be the blowup of of $\Rbar$ along $\Ibar$.
Let 
$M = M(I)$ as in 
Proposition-Definition~\ref{basic setup ell=1}.
 
We have:
\begin{enumerate}

 \item $M$ is a direct sum of maximal Cohen-Macaulay modules on the components of $\tilde R$, the module
 $M$ is $\omega_{\tilde R}$-self-dual, and $M$ has a linear resolution as an $S$-module.
 
 \item $\Rbar$ is equidimensional if and only if $M$ is a  maximal Cohen-Macaulay module over $\Rbar$.
 In this case
 $M$ is $\omega_{\Rbar}$-self-dual, and an Ulrich module over $\Rbar$, and
 $\dim R = \dim \Rbar$ if and only if $M$ is $\omega_R$-self-dual.
  
  \item 
  As a graded $R$-module, $\End(M) \cong M = I^\rho (\rho\delta) \cong \Ibar^{\rho}(\rho\delta)$ for $\rho\gg 0$. 
As a graded ring, $\End_R(M)\cong \tilde R$, and this ring is regular. In particular $\tilde R$ is the normalization of $\Rbar$, and the conductor of $\Rbar$  is $\Rbar \abar^{\rho}:\Ibar^{\rho}$ for $\rho\gg 0$. 

\end{enumerate}
\end{theorem}

\smallskip
Part (3) of the theorem includes the assertion that $M$ is a commutative algebra. The following lemma enables us to exploit this fact:

\begin{lemma}\label{lin pres ring}
Let $S$ be a standard graded polynomial ring over a field $k$. Suppose that $T$ is a non-negatively graded $S$-algebra that is finite,  generated in degree 0, and linearly presented as a module over $S$.
If $T_0$ is reduced, then $T$ is a product of standard graded polynomial rings over fields.
\end{lemma}

\begin{proof} The ring $T_0$ is finite over $k$ and hence 
Artinian. Thus
 $T_0$ decomposes as a product of fields $e_iT_0$, where the $e_i$ are orthogonal idempotents. Since it is a summand of $T$,
 the $S$-module $e_iT$ is also
 linearly presented, so it suffices to consider the case when $T_0 = K$ is a field. Let $\tilde S = K\otimes_kS$, and note that
 $T$ is naturally a cyclic $\tilde S$-module, generated in degree 0.
 
 We next show that $T$ is linearly presented as an $\tilde S$-module. Let $\tilde F$ be a 
 minimal homogeneous $\tilde S$-free resolution for $T$. Using an $S$-module splitting
 $\tilde S = S^m$ we see that $\tilde F$ is also an $S$-free resolution, which we call $F$.  Choosing bases in $\tilde F$ and corresponding bases in $F$,
 we see that the entries in the  matrix representations of the maps in $F$
 have the same positive degrees as occur in $\tilde F$, and thus $F$ is minimal. Since the initial map of $F$ is linear by hypothesis, the initial
 map of $\tilde F$ is linear.

Thus  $T$ is a cyclic, linearly presented $\tilde S$ module.
Therefore $T$ may be written as $\tilde S/m$, where $m$ is an ideal generated by linear forms, so $T$ is a polynomial ring.
\end{proof}

\vskip .1in
\begin{proof}[Proof of Theorem~\ref{stable residuals and linearly presented powers}]
We may assume that $I\not\subset \sqrt 0$, because otherwise $I^{\rho}$ and $\Rbar$ are zero and the assertions are trivial.
Now also $\Ibar$ is not nilpotent, $\ell(\Ibar)=\ell(I) =1$, and the ring $\Rbar$ is nonzero and reduced.

\noindent{\bf (3)} Note that $M$ is generated in degree 0 as an $S$-module and is linearly presented. Further, by Proposition-Definition~\ref{basic setup ell=1}, 
$M = \tilde R$ is an $S$-algebra. Since $\tilde R$ is reduced, Lemma~\ref{lin pres ring}  shows that
$\tilde R$ is regular. The rest of part (3) now follows from Proposition-Definition~\ref{basic setup ell=1}.

\smallskip

\noindent{\bf (1)} This follows from (3): We have $M \cong \tilde R$. Any regular ring
is Cohen-Macaulay, and, since $\tilde R$ is
Gorenstein, we have
$$
\tilde R = \Hom_{\tilde R}(\tilde R, \tilde R) \cong
\Hom_{\tilde R}(\tilde R, \omega_{\tilde R})
$$
up to a shift.

By Lemma~\ref{lin pres ring}, $\tilde R$ is a product of rings, each of which is obtained from $S$ by extending the ground field
and factoring out linear forms. Such a ring has a linear resolution as an $S$-module.
\smallskip

\noindent{\bf (2)} The ring $\tilde R$ is a finitely generated $R$-module by (3). The ring $\Rbar$ is equidimensional if and only if each component of
$\tilde R$ has maximal dimension as an $\Rbar$-module. By (2), this condition is equivalent to the condition that
$M$ is a maximal Cohen-Macaulay module over $\Rbar$. 

Moreover, 
if each  component of $\tilde R$ has maximal dimension as an $\Rbar$-module,
then $\omega _{\tilde R} \cong \Hom_{\Rbar}(\tilde R,\omega _{\Rbar})$. Thus
$$
\Hom_{\Rbar}(M, \omega_{\Rbar})
\cong
\Hom_{\tilde R}(M, \omega_{\tilde R}) \cong M,
$$
where the second isomorphism is part of  (1). If $\dim R = \dim \Rbar$, these statements
follow for $R$ in place of $\Rbar$, and the converse is obvious.

The fact that $M$ is Ulrich follows from the linearity statement in (1).
\end{proof}

The following lemma shows that the linearity hypothesis in Theorem~\ref{main} is preserved in the reduction
to Theorem~\ref{stable residuals and linearly presented powers}.
\begin{lemma}\label{lin pres modules}
 Let $S$ be a standard graded polynomial ring over an infinite field,  
 and let $I\subset S$ be a homogeneous ideal, minimally generated by forms of degree $\delta$.
 Let $a_1,\dots, a_{n}$ be homogeneous elements of $I$ of  degree $\leq \delta+1$. 
 For  $0\leq \nu\leq n$, set $J_\nu = (a_1,\dots, a_\nu)\subset S$. 
 
 If the high powers of $I$
 have linear presentation, then for all $0\leq \nu\leq n$ and all sufficiently large $\rho$ the module
 $I^\rho/J_\nu I^{\rho-1}$ has linear presentation as an $S$-module.
\end{lemma}

\begin{remark}
 If the $a_{i}$ are chosen sequentially generally, the same proof shows that if high powers of $I$ have linear resolution then for all $0\leq \nu\leq n$ and all sufficiently large $\rho$ the module
 $I^\rho/J_\nu I^{\rho-1}$ has linear resolution as an $S$-module.
\end{remark}

\begin{proof}
We do induction on $\nu$, the case $\nu=0$ being the hypothesis. We assume the result for $\nu$.
For  $\gamma\geq 1$ and $0\leq \nu\leq n-1$  the sequences
 $$
 (I^{\gamma-1}/J_{\nu} I^{\gamma-2})(-\delta') \rTo^{a_{\nu+1}} I^{\gamma}/J_{\nu} I^{\gamma-1} \to I^{\gamma}/J_{\nu+1} I^{\gamma-1} \to 0
 $$
are exact, where $\delta\leq \delta' = \deg a_{\nu+1}\leq \delta+1$. Using the minimal homogeneous  generators of the middle term $I^{\gamma}/J_{\nu} I^{\gamma-1}$, we get a possibly non-minimal set of generators of 
$I^{\gamma}/J_{\nu+1} I^{\gamma-1}$ having degree $\gamma\delta$, with relations of degrees 
$\gamma \delta+1$ and $(\gamma-1)\delta +\delta'\leq \gamma\delta+1$.
This implies that 
$I^{\gamma}/J_{\nu+1} I^{\gamma-1}$
has linear presentation.
\end{proof}

\begin{proof}[Proof of Theorem~\ref{main}] Write $d=\dim S$ and let $a_1, \ldots, a_{\ell -1}$ be the $\ell -1$ general elements
of $I$ that generate $J$. We apply Theorem~\ref{stable residuals and linearly presented powers} to the ideal
$IR\subset R$, and verify the hypotheses as follows:

\begin{itemize}
 \item $\ell(IR) \leq 1$: This holds because the general elements $a_1, \ldots, a_{\ell -1}$ are part of a minimal generating set for a minimal reduction of $I$.
  
 \item $M$ is linearly presented as an $S$-module: By Lemma~\ref{basic setup higher ell},
$M \cong (I^{\rho}/JI^{\rho-1})(\rho\delta)$ for large $\rho$. By Lemma~\ref{lin pres modules},  this
is linearly presented as an $S$-module.
 \end{itemize}
To apply item (2) of Theorem~\ref{stable residuals and linearly presented powers}, 
we will show that $\Rbar$ is equidimensional of dimension $d -\ell +1$: 

We first argue that $\Rbar\neq0$. Indeed, if $\Rbar=0$ then 
 $IR \subset \sqrt 0_R$, hence $I^n \subset J:I$ for some $n >0$ and therefore $I^{n+1} \subset J$. 
 As $I \subset \sqrt J$ and $J$ is generated by $\ell -1$ general elements of $I$, it follows from \cite[Theorem 4]{S} that
 $\ell -1 \geq \ell$, a contradiction.  
 
 Next, any minimal prime of the ring $\Rbar$ arises from a minimal prime $Q$ of the ideal
 $K:=J:I$ that does not contain $I$, and we need to show that $\codim Q = \ell -1$ or, equivalently, $\codim K_Q=\ell -1$. But
 $K_Q=J_Q$ since $I \not\subset Q$. Finally, the generators $\frac{a_1}{1}, \ldots, \frac{a_{\ell-1}}{1}$ of $J_Q$ form a 
 regular sequence on $S_Q$, because
 the general elements $a_1, \ldots, a_{\ell-1}$ of $I$ are a filter regular sequence with respect to $I$ and $Q \not\supset I$. 
 
 Thus $\codim K_Q= \codim J_Q =\ell -1$. Together with item (2) of Theorem~\ref{stable residuals and linearly presented powers}, this gives the assertion
 of item {\bf (1)}.

We now prove item {\bf (2)}. We have seen in the previous step that
 $\dim \Rbar=d-\ell +1\leq \dim R$. Thus $\dim \Rbar =\dim R$ if and only if $\codim K\geq\ell -1$.  
 
 Finally, item (3) of Theorem~\ref{stable residuals and linearly presented powers} now implies {\bf (3)}.
\end{proof}
 
\vspace{.05in}
\section{The minors of a  $2\times n$ matrix of linear forms}

Throughout this section we assume $n\geq 2$. We say that an ideal $I \subset S$ is \emph{$s$-residually $S_2$} if,
for every $i\leq s$ and every $i$-residual intersection $K$ of $I$, the ring $S/K$ satisfies Serre's condition $S_2$. See \cite{CEU} for more information.

 \begin{theorem}\label{2xn prelim}
Let $k$ be a field of characteristic 0. Suppose that $I \subset S= k[x_{1,1},\dots,x_{2,n}]$ is the ideal of $2\times 2$ minors of the generic matrix 
 $$
\begin{pmatrix}
 x_{1,1}&\dots&x_{1,n}\\x_{2,1}&\dots&x_{2,n}
\end{pmatrix} .
$$
Let $\ell := \ell(I)$, which is equal to $2n-3$ by by Proposition~\ref{matrix facts}.
The ideal $I$ is $(\ell-2)$-residually $S_{2}$. 
In particular, if $s \leq \ell-1$ and $K$ is an s-residual intersection of $I$, then
$K$ is unmixed of codimension exactly $s$. If, in addition, the residual intersection
is geometric, then the image of $I$ in $S/K$ contains a non-zerodivisor.
\end{theorem}

\begin{proof}
Note that $I$ is a complete intersection on the punctured spectrum of $S$ and thus satisfies $G_{2n}$. 
By~\cite[Theorem 4.3]{RWW}, $\Ext^{n+j-1}_{S}(S/I^{j}, S) =0$ for $1\leq j\leq n-3 = (\ell-2)-\codim I +1$ (this is where we require characteristic 0). By~\cite[Corollary 4.2]{CEU}, this implies that $I$ is $(\ell-2)$-residually $S_{2}$. 

 From~\cite[Proposition 3.1]{CEU} we know that $I$ is $(\ell-1)$-parsimonious. Note that $K$ is a proper ideal because $s$ is less than the number of generators of $I$. Thus we may apply
 \cite[Proposition 3.3(a)]{CEU} and conclude that $K$ is unmixed of codimension exactly $s$. If, in addition, $K$ is a
geometric residual intersection, then $\codim I+K\geq s$ so $I$ is not in any associated prime of $K$. 
\end{proof}

\begin{corollary}\label{unmixedness}
 Suppose that $I$ is the ideal of $2\times 2$ minors of a $2\times n$ matrix $A$ over a local ring $T$ containing a field
 of characteristic 0, and 
 assume that $\codim I = n-1$. If $s \leq 2n-4$ and $K$ is an s-residual intersection of $I$, then
 every minimal prime of $K$ has codimension exactly $s$. If, in addition, the residual intersection
is geometric, then  $I$ is in no minimal prime of $K$.
 \end{corollary}

\def \tA{{\tilde A}}
\def\tI{{\tilde I}}
\def\tJ{{\tilde J}}
\def\tK{{\tilde K}}
\def\tT{{\tilde T}}
\def\lbracket{{[\kern-1.5pt[}}
\def\rbracket{{]\kern-1.5pt]}}

\begin{proof} We may assume that the entries of $A$ are in the maximal ideal of $T$, since
otherwise $I$ is a complete intersection.

Let $J \subset I$ be an ideal with $s$ generators such that $K = J:I$.
Let $\tT = T\lbracket x_{1,1},\dots,x_{2,n}\rbracket$, and let $\pi:\tT \to T$ be a map sending $x_{i,j}$ to the $(i,j)$ entry $A_{i,j}$ of $A$, and note that the kernel of $\pi$ is generated by the regular sequence $\alpha_{i,j} = x_{i,j}-A_{i,j}$.
Let $\tI$ be the ideal of $2\times 2$ minors of the generic $2\times n$ matrix 
 $$
\begin{pmatrix}
 x_{1,1}&\dots&x_{1,n}\\x_{2,1}&\dots&x_{2,n}
\end{pmatrix},
$$
 so that $\pi(\tI) = I$. Let 
$\tJ\subset \tI$ be an ideal with $s$ generators such that $\pi(\tJ) = J$, and let $\tK = \tJ:\tI$.

Since $\tT/\tI$ is Cohen-Macaulay, and the codimension of $\tI$ is equal to that of $I$, the $2n$ elements
$\alpha_{1,1},\dots,\alpha_{2,n}$ form a regular sequence on $\tT/\tI$. It now follows from 
\cite[4.1]{HU1} that $\sqrt{\pi(\tK)} = \sqrt K$. Thus
$$
\codim \tK \geq \codim \pi(\tK) = \codim K \geq s,
$$
so $\tK$ is an $s$-residual intersection of $\tI$. As $s\leq 2n-4$, Theorem~\ref{2xn prelim} implies that $\tK$ is 
unmixed of codimension exactly $s$.  

Since $\codim \pi(\tK)\geq s$, it follows that the sequence $\alpha_{1,1},\dots,\alpha_{2,n}$ is part of
a system of parameters of $\tT/\tK$, and thus all minimal primes of $\pi(\tilde K)$ have codimension exactly $s$.

Using $\sqrt{\pi(\tK)} = \sqrt K$ again, we see that all minimal primes of $K$ have codimension exactly $s$.

The last statement follows immediately.
\end{proof}

\begin{theorem}\label{scrolls}
Suppose that $I$ is the ideal of $2 \times 2$ minors of a $2\times n$ matrix $A$ of linear forms in a polynomial ring $S$
over a field of characteristic 0, and suppose that the entries of $A$ span a vector space of dimension $c$. 

If $I$ has codimension $\min\{n-1, c\}$, then the hypotheses and conclusions of 
Theorem~$\ref{main}$ hold for $I$, and the ring $\Rbar$ in Theorem~$\ref{main}$ is $R_{\rm red}=R/\sqrt 0$.
\end{theorem}

\begin{proof} Theorem~\ref{list}(1) and Remark~\ref{Bertini} imply that the hypotheses and conclusions of Theorem~\ref{main}(\ref{2.1.1}) and (\ref{2.1.3}) hold for $I$.

The ideal $I$ satisfies $G_{c}$, and $\ell:= \ell(I) \leq c$. Thus $K$ is a geometric $\ell-1$ residual intersection. Hence we may apply Theorem~\ref{main}(\ref{2.1.2}).

If $\codim I = n-1$, then Corollary~\ref{unmixedness} shows that $I$ is not contained in any minimal prime of $K$. On the other hand, if $\codim I = c$ then $c=\ell$, so $K$ is a complete intersection of codimension $c-1$, and again $I$ is not contained in any minimal prime of $K$.

Since $R$ is reduced away from $V(I)$, in follows in both cases that 
$$
\Rbar = R/(0:I^\infty) = R/(\sqrt 0:I^\infty) = R/\sqrt 0.
$$
 \end{proof}

In the case of the generic $2\times n$ matrix, we can be very explicit. 
\def\tR{{\tilde R}}

  \begin{theorem}\label{2xn} Let $S= k[x_{1},\dots,x_{n}, y_{1},\dots,y_{n}]$ be a polynomial ring over
  a field $k$ of characteristic 0. Suppose that $I \subset S$  is the ideal of $2\times 2$ minors of the generic matrix
  $$
\begin{pmatrix}
 x_{1}&\dots&x_{n}\\y_{1}&\dots&y_{n}
\end{pmatrix}.
$$
 The ideal $I$ has analytic spread $\ell := \ell(I) = 2n-3$ and reduction number $r := r(I)= n-3$ by Proposition~\ref{matrix facts}.
 
Let $J$ be an ideal
  generated by $\ell-1$ general quadrics $a_{1}, \dots a_{\ell-1}$ in $I$. Set $R = S/(J:I)$, $\Rbar = R/(0:I^\infty)$, and $M = M(IR)$. In addition to the assertions
  of Theorem~\ref{main} we have:
 \begin{enumerate}

\item $R$ has an isolated singularity, and $\Rbar = R$.

\item If $a$ is a general quadric in $I$ and $\rho \geq r$ then $a^{\rho-r}: I^rR(2r) \to I^\rho R(2\rho)$ is an isomorphism, so $M = I^rR(2r)$.

\item $M\cong \omega_{R}(4).$

\item If $k$ is algebraically closed, then $J:I$ is the intersection of $\frac{1}{n-1}{2(n-2)\choose n-2}$
linear prime ideals of codimension $2n-4$. These may be described as follows: the generators
$a_{1}, \dots, a_{2n-4}$ may be regarded as linear forms in $k[I_{2}]$, which may be identified with the homogeneous
coordinate ring of the Grassmanian $G(2,n)$ of $n-2$-dimensional subspaces of $k^{n} = \langle x_{1}\dots x_{n}\rangle$. 
Since the $a_{i}$ are general, the plane cut out by these forms intersects the Grassmanian in $\frac{1}{n-1}{2(n-2)\choose n-2}$ reduced points.
Each of these points corresponds to a subspace $\langle L_{1}(x),\dots, L_{n-2}(x)\rangle$, which yields a linear prime ideal
$$
(L_{1}(x),\dots, L_{n-2}(x), L_{1}(y),\dots, L_{n-2}(y))
$$
that is a minimal prime of $J$.
\end{enumerate}
  \end{theorem}

\begin{proof}
 For $i\leq \ell-1$ we set $J_{i}= (a_{1},\dots,a_{i})$. Since $I$ satisfies $G_{2n}$, the ideal $J_i:I$ is a geometric residual intersection of $I$ and is unmixed of codimension  $i$ by Theorem~\ref{2xn prelim}.

\noindent{\bf (1)}
 By Theorem~\ref{2xn prelim} the image of $I$ in $R$ contains a nonzerodivisor. Thus 
 $R$ is reduced by Remark~\ref{Bertini}, and therefore $R = \Rbar$ by Theorem~\ref{scrolls}.

By \cite{HU2}, the generic residual intersection
is nonsingular in codimension 3, and since general fibers are reduced, Bertini's theorem implies
that general fibers are nonsingular in codimension 3. 

Since $R$ has dimension 4,  its singularity is isolated. 

\noindent{\bf (2)} Let $a\in I$ be a general quadric. By Theorem~\ref{2xn prelim} the element $a$ is a nonzerodivisor on $R$. The ideal $IR$ has analytic spread at most 1 and reduction number at most $r$. Thus if $\rho\geq r = r(I)$, then $a^{\rho-r}: I^rR(2r) \to I^\rho R(2\rho)$ is an isomorphism, so $M = I^rR(2r)$.

\noindent{\bf (3)}  Let $M_{j,i} := I^{j}/J_iI^{j-1}(2j)$, which is generated in degree 0. We will show that $M \cong M_{r,\ell-1}\cong M_{\rho,\ell-1}$ for all $\rho\geq r$. Moreover, we will show that $M_{r+1,\ell-1} \cong \omega_R(4)$. For this we must estimate the depth of $ M_{r,\ell-1}$, and for this in turn we first prove that certain syzygies of  $M_{j,i}$ have low degree.

\begin{lemma}\label{linearity} 
If  $\, 1\leq j\leq r$ and $\, 0\leq i\leq n+j-1$ then 
$\depth M_{j,i} \geq \min(4,n)$, and the $m$-th free module in a minimal graded $S$-free resolution of $M_{j,i}$ is generated in degree $m$ for all $m>i$.
\end{lemma}

\begin{proof} We adopt the notation of the proof above.
We first consider the case $j=1$. Set $g = \codim I = n-1$. We must treat the cases with $i\leq n = g+1$. 
For $i<g$ we have $J_{i}:I = J_{i}$, a complete intersection of codimension $i$. If $i = g$ then, since $S/I$ is Cohen-Macaulay, the link $S/(J_{g}:I)$ is also Cohen-Macaulay. By~\cite[Lemma 1.7]{U},  then the depth of $S/J_{g+1}$ is at least $\dim S -g -1 = n$. In each of these  cases the length of the resolution of $M_{1,i} = (I/J_{i})(2)$ is at most $i$, proving both statements.

We now do induction on $i$.

If $i=0$, then $M_{i,j} = I^{j}(2j)$. By \cite[Theorem 5.4 and the beginning of its proof]{ABW} the minimal graded free resolution 
of $I^{j}$ is linear and of length at most $2n-4$ for every $j$, as required.

We now suppose $i>0$. We may assume that $j>1$. Consider the sequence
$$
0\to M_{j-1,i-1} \rTo^{\alpha} M_{j,i-1} \rTo M_{j,i} \to 0,
$$
where $\alpha$ is multiplication by $a_{i}$. It follows from the definitions that the sequence is right exact. We will show that it is exact.

Since $I$ is a complete intersection on the punctured spectrum, \cite[Lemma 2.7(a)]{U} with $s = 2n$ shows that the left-hand map in this sequence
is a monomorphism locally on the punctured spectrum because $r\leq 2n-(n-1)+2$. By induction $M_{j-1,i-1}$ has positive depth, so the
sequence in also left exact as claimed. 

\def\bH{{\bf H}}
Let $\tilde \alpha: F_{\bullet} \to G_{_{\bullet}}$ be the map of minimal graded free resolutions induced by $\alpha$. The
minimal graded free resolution $H_{_{\bullet}}$ of $M_{j,i}$ is a direct summand of the mapping cone of 
$\tilde \alpha$. Hence it follows by induction that $H_{m}$ is generated in degree $m$ for all $m>i$.

Finally, we must show that the length of $H_{\bullet}$ is at most $2n-4$. By the induction hypothesis, the length is at most $2n-3$. Further,
$H_{2n-3} $ is a direct summand of $F_{2n-4}$. Moreover, $F_{2n-4}$ is generated in degree $2n-4$
because $2n-4 =n+r-1\geq n+j-1>i-1$. Thus $H_{2n-3}$ is generated in degree $2n-4$.
On the other hand,
$H_{2n-4}$ is generated in degrees at least $2n-4$, so the  map $H_{2n-3}\to H_{2n-4}$ in $H_{\bullet}$
is given by a scalar matrix. Since $H_{\bullet}$ is minimal by assumption, it follows
that $H_{2n-4}$ is in fact 0, as required.
\end{proof}

Continuing with the proof of part (3), we have
a natural surjection of $R$-modules $\pi: M_{r, \ell-1} = (I^{r}/JI^{r-1})(2r) \to (IR)^{r}(2r) = M$. Recall that $J:I$ is a geometric $(\ell-1)$-residual intersection. Moreover, on the punctured spectrum $I$ is a complete intersection, hence by \cite[Lemma 2.7(c)]{U} with $s = 2n$, the kernel of $\pi$ is 0 on the punctured spectrum. On the other hand, $M_{r, \ell-1}$ has depth $\geq 1$ by Lemma~\ref{linearity}, so the kernel is 0 and we see that
 $\pi$ is an isomorphism.

Let $a\in I$ be a general quadric, and consider the diagram: 
$$
\begin{diagram}
M_{r,\ell-1}&\rTo^{\cong}_\pi&(IR)^{r}(2r)\\
\dTo^{a^{\rho-r}}&&\dTo_{a^{\rho-r}}\\
M_{\rho,\ell-1}&\rTo& (IR)^{\rho}(2\rho)
\end{diagram}
$$
By item (2) the right hand vertical map is an isomorphism. It follows that the left-hand vertical map is a monomorphism.
 For $\rho \geq r = r(I)$ we have 
$$
I^\rho = (J,a)^{\rho -r}I^r = JI^{\rho-1}+a^{\rho-r}I^r,
$$
 so the left hand vertical map is also a surjection. Thus all the maps in the square are
 isomorphisms, so $M = M_{r,\ell-1} = M_{r+1,\ell-1}$.
 \smallbreak

By Theorem~\ref{2xn prelim} the ideal $I$ is $(\ell-2)$-residually $S_2$, so by~\cite[Proposition 3.1]{CEU} it is ($\ell-1$)-thrifty.
In particular it is $(\ell-2)$-thrifty, and thus by~\cite[Lemma 2.4(b)]{CEU},
$$
(J_{i-1}:I, a_i):I = J_i:I
$$
for all $i\leq \ell-1$. By Theorem~\ref{2xn prelim} the right hand side, and thus also the left hand side, has codimension exactly $i$. This verifies the hypothesis of~\cite[Theorem 4.1]{EU}, and thus 
there is a natural homogeneous map
 $$
 \mu: \bigl(I^{(\ell-1)-g+1}/J_{\ell-1}I^{(\ell-1)-g}\bigr)(2(\ell-1)-2n) \rTo \omega_{S/(J_{\ell-1}:I)} = \omega_R
 $$ 
that is an isomorphism on the punctured spectrum since $I$ is locally a complete intersection there. (Though~\cite[Theorem 4.1]{EU} was proven in the local case, the twists can be recovered from the proof.)

We have $\ell-1-g+1 = n-2 = r+1$, and $2(\ell-1)-2n = 2(r+1) - 4$ so
the source of $\mu$ is $M_{r+1, \ell-1}(-4)$. Since this module has depth $\geq 2$ by
Lemma~\ref{linearity}, it follows that $\mu$ is an isomorphism, proving (3).
 
\smallbreak
\noindent{\bf (4)} 
We next prove that the linear prime ideals described in (4)  contain $J:I$. A point $z$ on the Grassmannian corresponds
to a $2\times n$ matrix of rank 2, which, after coordinate transformation, may be taken to be 
$$
\begin{pmatrix}
 0&\dots&0&1&0\\
  0&\dots&0&0&1\\
\end{pmatrix}.
$$
The Pl\"ucker coordinates of $z$ are then 
$$
p_{\mu,\nu} = \begin{cases}
1 \quad\hbox { if } (\mu,\nu) = (n-1,n)\,\\
0 \quad\hbox{ otherwise .}\\
\end{cases}
$$
We may write the $a_{i}$ in the form $\sum \lambda_{\mu,\nu}^{i}p_{\mu,\nu}$. To say the point $z$ is on the linear section defined by the $a_{i}$ means that the coefficients $\lambda_{n-1,n}^{i}$ are all 0. Thus the
$a_{i}$ are in the ideal 
$
L := (x_{1},\dots,x_{n-2}, y_{1},\dots,y_{n-2})
$, which is the prime corresponding to $z$.
Finally, this implies that $J:I\subset L$ because $I\not\subset L$.
In particular this shows that the multiplicity of $R$ is at least the degree of the Grassmanian, which is
$\frac{1}{n-1}{2(n-2)\choose n-2}$.  

The degree $2r$ component of the graded module $(IR)^{r}$ is a homomorphic image of the degree $r$ component of $k[I_{2}]/((a_{1}, \dots,  a_{2n-4})$. 
The $a_{i}$ are general linear forms in $k[I_{2}]$, the coordiinate ring of the Grassmannian in the Pl\"ucker embedding. Because this ring is Cohen-Macaulay of dimension $2n-3$, the ring $k[I_{2}]/(a_{1}, \dots,  a_{2n-4})$ is a one-dimensional Cohen-Macaulay ring of multiplicity equal to the degree of the Grassmannian, and thus the number of generators of $(IR)^{r}$ is bounded by the degree of the Grassmannian.

We deduce
that the multiplicity of $R$ is equal to the number of linear minimal primes as above. Since $R$ is unmixed
of codimension $2n-4$,
this shows that $J:I$ is the intersection of these linear primes proving (4).
\end{proof}

\section{Determinantal ideals}

Theorem 3.7 of~\cite[Theorem 3.7]{BCV} gives a large family of determinental ideals whose powers have linear resolutions, reproduced in part (2) of the following:

 \begin{theorem}\label{Bruns et al}
Suppose that $A$ is an $m\times n$
matrix of linear forms in a polynomial ring $S$, with $m\leq n$, and suppose that the entries of $A$ generate a vector space of dimension $c$.
Let $I$ be the ideal of $m\times m$ minors of $A$. If either:
\begin{enumerate}
 \item $\codim I = c$; \ or 
 \item  $\codim I = m-n+1$ and for  $2\leq k\leq m-1$ the ideal of $k\times k$ minors of $A$ has codimension
 $\geq \min((m-k+1)(n-m)+1, c)$;
\end{enumerate}
then every power of $I$ has a linear resolution.
\end{theorem}

Since the powers of the ideal of the Veronese surface also have linear resolutions (\cite[Proposition 3.12]{BCV}) the powers of the ideal of every variety of minimal degree have linear resolutions. 

It seems plausible that, in the setting above, if $I$ itself has linear presentation (respectively, linear resolution), then all its powers do too.  In the case $m=2$, the condition for
$I$ itself to have linear presentation or resolution is known in terms of the Kronecker classification of linear $2\times n$ matrices; see \cite{C-J}. In fact, the condition that high appears to be more general still: for example, let $I$ be ideal of 
$2\times 2$ minors of the matrix
$$
\begin{pmatrix}
 0&x_{1}&\cdots&x_{5}&|&y_{0}&y_{1}&y_{2}\\
x_{1}&\cdots&x_{5}&0&|&y_{1}&y_{2}&y_{3}\\
\end{pmatrix}.
$$
\smallbreak
According to Macaulay2\cite{M2}, the Betti tables of the first 3 powers of $I$ (in characteristic 101) are:
\smallbreak
\tiny{
\tt{
\begin{verbatim}
o4 = total: 1 36 169 386 531 470 271 99 21 2
         0: 1  .   .   .   .   .   .  .  . .
         1: . 36 169 383 514 430 221 64  8 .
         2: .  .   .   3  17  40  50 35 13 2

            0   1    2    3     4     5    6    7   8   9
o5 = total: 1 414 2542 7124 11754 12395 8514 3708 934 104
         0: 1   .    .    .     .     .    .    .   .   .
         1: .   .    .    .     .     .    .    .   .   .
         2: .   .    .    .     .     .    .    .   .   .
         3: . 414 2542 7124 11752 12385 8494 3688 924 102
         4: .   .    .    .     2    10   20   20  10   2

            0    1     2     3     4     5     6     7    8   9
o6 = total: 1 2544 17028 50967 88676 97776 69804 31458 8172 936
         0: 1    .     .     .     .     .     .     .    .   .
         1: .    .     .     .     .     .     .     .    .   .
         2: .    .     .     .     .     .     .     .    .   .
         3: .    .     .     .     .     .     .     .    .   .
         4: .    .     .     .     .     .     .     .    .   .
         5: . 2544 17028 50967 88676 97776 69804 31458 8172 936

\end{verbatim}
}}

\normalsize
\noindent and the 4th power also has linear resolution, suggesting that higher powers will too.
\begin{proof}
Suppose first that $\codim I = c$, so that in particular $c<= m-n+1$. We may harmlessly assume that the entries of $A$ span the space of all linear forms and that the ground field is infinite.
We may write $S$ as $T/J$ where $T$ is a polynomial ring in $mn$ variables in such a way that $A$ is the specialization of a generic matrix $B$. For a generic choice of intermediate specialization $T'$ of dimension $m-n+1$, with
$$
T\twoheadrightarrow T' \twoheadrightarrow S
$$
the ideal of $m\times m$ minors $I'$ of the specialization $B'$ of $B$ to $T'$ will have codimension $n-m+1$. It follows that the minimal resolution of $I'$ is the Eagon-Northcott complex, and thus the ${n\choose m}$ minors of $B'$ are linearly independent. Since the vector space dimension of the degree $m$ component of $T'$ is also ${n\choose m}$,
the ideal $I'$ is the $m$-th power of the maximal ideal of $T'$. Specializing further to $S$ we see that $I$ is the $m$-th power of the maximal ideal. 

The sufficiency of $(2)$ is the main theorem of Bruns, Conca and Varbaro~\cite{BCV}. 
\end{proof}

\subsection*{Generic matrices}

The analytic spread and reduction number of an ideal of maximal minors of a generic matrix are known; for the reader's convenience we reproduce the result.

\begin{proposition}\label{matrix facts}
Let ${\mathcal X}$ be the generic $m\times n$ matrix of variables over the ring $S:= k[x_{1,1},\dots, x_{m,n}]$, with $m\leq n$, and let $I = I_{m}({\mathcal X})$ be the ideal of $m\times m$ minors. The analytic spread
 of $I$ is $\ell(I) = m(n-m)+1$ and, when the ground field $k$ is infinite, the reduction number of $I$ is $r(I) = \ell(I) -n$.
\end{proposition}
\begin{proof}
 Let $\mm\subset S$ be the ideal generated by the entries $x_{i,j}$ of $\mathcal X$. 
 The special fiber ring ${\mathcal F(I)} := S/\mm \oplus I/\mm I \oplus I^2/\mm I^2\cdots$ of $I$ is the homogeneous coordinate ring of the Grassmannian $G(m,n)$ in its Pl\"ucker embedding. Since $G(m,n)$ is a variety of dimension $m(n-m)$, the analytic spread of $I$ is $\ell(I) = m(n-m)+1$.

Now assume that the ground field is infinite.  The reduction number $r(I)$ of ${\mathcal F(I)}$ is the maximal degree of a socle element after reducing ${\mathcal F(I)}$ modulo a general linear system of parameters \cite[Theorem 8.6.6]{SH}. Because the homogeneous coordinate ring of the Grassmannian is Cohen-Macaulay, we can relate this to the degree of the generators of the canonical module. The canonical module of the Grassmannian $G(m,n)$ is 
$\mathcal {O}_{G}(-n)$ in the Pl\"ucker embedding (see for example \cite[Proposition 5.25]{EH}). Thus modulo a general sequence of $\ell = m(n-m)+1$ linear forms, the socle is in degree
$\ell(I)-n$, and the reduction number is thus $r(I) = \ell(I) -n$.
\end{proof}

It is interesting to ask when an ideal of maximal minors has an $(\ell-1)$-residual intersection, so that part (2) of Theorem~\ref{main} is non-trivial. We thank Monte Cooper and Edward Price for pointing out an error in a previous version.

\begin{proposition} Let ${\mathcal X}$ be the generic $m\times n$ matrix of variables over the ring $k[x_{1,1},\dots, x_{m,n}]$, with $m\leq n$, and let $I = I_{m}({\mathcal X})$ be the ideal of $m\times m$ minors. Let $\ell :=\ell(I)$, which is $m(n-m)+1$ by Proposition~\ref{matrix facts}.
\begin{enumerate}
 \item The ideal $I$ satisifies $G_{\ell}$  if and only if one of the following holds: 
\begin{itemize}
 \item $m\leq 2$;
 \item $n\leq m+2$;
 \item $n = m+3$ and $m\leq 5$.
\end{itemize}
\item The ideal $I$ satisifies $G_{\ell-1}$ if and only if it satisfies $G_{\ell}$ or
\begin{itemize}
 \item $m=3,\ n=7$;
  \item $n = m+3$ and $m\leq 6$.
 \end{itemize}

\item $I$ does not have any $(\ell-1)$-residual intersection if one of the following holds:
\begin{itemize}
 \item $n = m+3$ and $m= 10$ or $11$ or $m\geq 14$;
 \item $n= m+4$ and $m\geq 6$;
 \item $n\geq m+5$ and $m\geq 3$.
\end{itemize}
\end{enumerate}
 \end{proposition}
\begin{proof} For every prime $P \in V(I)$ one has $P\in V(I_{t+1}(A))\setminus V(I_{t}(A))$ for some $t$ 
with $0 \leq t \leq m-1$, and the  minimal number of generators of $I_{P}$ is exactly 
${n-t\choose m-t}$.

Thus  the condition $G_{s}$  holds for $I$ if and only if
$$
{n-t\choose m-t} \leq \codim I_{t+1}(A) = (m-t)(n-t) 
$$
whenever $\codim I_{t+1}(A)\leq s-1 $. Given this, the verification of items (1) and (2) is not difficult.

If $I$ admits an $(\ell -1)$-residual intersection, then locally in codimension $\ell -2$, the ideal $I$ can be generated
by $\ell -1$ elements. In other words, 
$$ 
{n-t\choose m-t} \leq \ell -1 = m(n-m)
$$ 
whenever  $\codim I_{t+1}(A) = (m-t)(n-t) \leq \ell -2$. Again, part (3) follows easily from this.
\end{proof}

\vspace{.05in}

\section{Implications and special cases of the conjectures}

\subsection{Implications of Conjecture~\ref{main conjecture0}\label{proof of 0.5}}
\begin{proof}[Proof of Proposition~\ref{cor of conjecture}.]
We may assume that $k$ is infinite.

\noindent (1) The result is trivial if $I$ is not a complete intersection, so we assume that it is not. In this case, $\ell >c$ by \cite{CN}. Thus 
$\ell=c+1$. It follows that the ideal $J:I$ of Conjecture \ref{main conjecture0} is a link, hence unmixed, and
the ideal $IR$ is principal.  As $I$ is generically a complete intersection, the link is geometric and 
$IR$ is generated by a single non-zerodivisor. If $I^{\rho}R$ were a maximal Cohen-Macaulay $R$-module
for some $\rho >0$, then $R=S/(J:I)$ is Cohen-Macaulay, hence so is $S/I$ because the unmixed ideal $I$ is also a link of $J:I$. 

\noindent (2) Because $I$ is of linear type, $\ell$ is the number of generators of $I$. Let $\phi$ be a homogeneous presentation matrix of $I$ with respect to a general choice of homogeneous generators $f_1,\dots,f_\ell$ of $I$. The ideal $P$ defining the symmetric algebra of $I$ as a quotient of
$S':=k[T_{1},\dots, T_{\ell}]\otimes_{k}S$ is generated by the entries of the row vector $(T_{1}, \dots, T_{\ell}) \circ \phi$. 

Let $S'' = k(T_{1},\dots, T_{\ell})\otimes_{k}S$. Over $S''$, the row vector $(T_{1}, \dots, T_{\ell}) \circ \phi$ is the last row
of a presentation matrix of $IS''$ with respect to some homogeneous generators $g_1,\dots,g_\ell$.
Thus $PS''$ has the form $(g_1,\dots,g_{\ell-1}):IS''$. Since the $f_i$ were chosen generally over $k$, they are general over
$k(T_{1},\dots, T_{\ell})$, and it follows that the $g_i$ are general over $k(T_{1},\dots, T_{\ell})$. 

By hypothesis, $\Sym(I) = \RR(I)$, a domain of dimension $d+1$. Thus $PS''$ is a geometric $(\ell-1)$-residual intersection of $IS''$, and 
$I(S''/PS'')$ is generated by a nonzerodivisor. By  Conjecture~\ref{main conjecture0}, $I^{\rho}(S''/PS'')$ is a maximal Cohen-Macaulay module over $S''/PS''$. Since this is a principal ideal generated by a non-zerodivisor, 
$S''/PS'' =\Sym(I) \otimes_{S'} S'' =\RR(I)\otimes_{S'} S''$ is Cohen-Macaulay, and it follows that 
$\RR(I)_{(x_{1},\dots,x_{d})\RR(I)}$ is too.
\end{proof}

\subsection{Special cases of the conjectures}

\smallskip

The next result has been proven with an additional hypothesis in \cite[Theorem 2.6]{H1}.


\begin{theorem}\label{ACI}
Let $S$ be a local Gorenstein ring and let $I\subset S$ be an almost complete intersection ideal such that
$S/I$ is equidimensional. If 
$\depth \, (S/I)_P \geq \frac{1}{2} \,  \dim \,  (S/I)_P \, $ for every $P \in V(I) ,$ then $S/I$ is Cohen-Macaulay.
\end{theorem}

\begin{proof}
Let $J\subset I$ be a complete intersection of the same codimension as $I$ such that $I/J$ is cyclic, and consider $K = J:I$. 
Our assumptions imply that $I$ is unmixed. Therefore $I= J : K$ and it suffices to prove the Cohen-Macaulayness of $S/K$.

Notice that $\omega_{S/K}\cong I/J \cong S/K.$ Thus by \cite[Theorem 1.6]{HO} or \cite[Lemma 5.8]{H2} it suffices to show that
$$
\depth \, (S/K)_P \geq  1 + \tiny{\frac{1}{2}} \, \dim \, (S/K)_P 
$$
for every $P \in V(K)$ with $\dim \, (S/K)_P \geq 2 .$
We may assume that $P \in V(I)$ since otherwise $(S/K)_P=(S/J)_P$ is Cohen-Macaulay. But then
$\depth \, (S/I)_P \geq \frac{1}{2} \,  \dim \,  (S/I)_P$ and $\dim \,  (S/I)_P=\dim \,  (S/K)_P$. Now the 
exact sequence 
$$ 
0 \rightarrow S/K \cong I/J \longrightarrow S/J \longrightarrow S/I \rightarrow 0
$$
shows that $\depth \, (S/K)_P \geq  1 + \frac{1}{2} \,  \dim \,  (S/K)_P ,$ as required.
\end{proof}

\smallskip

Notice that an almost complete intersection ideal $I \subset S$ satisfies the assumptions of Theorem \ref{ACI} if $I$ is unmixed
and $S/I$ is almost Cohen-Macaulay, which means that $\depth S/I \geq \dim S/I - 1 .$
\begin{corollary}
If $I\subset S:= k[x_1,\dots, x_d]$ is an umixed monomial almost complete intersection, then $S/I$ is Cohen-Macaulay.
\end{corollary}
\begin{proof}
The Taylor resolution shows that the projective dimension of the polynomial ring modulo a monomial ideal is bounded by the number of generators of the ideal; thus any monomial almost complete intersection is almost Cohen-Macaulay.
\end{proof}

\begin{corollary}
 With hypotheses as in Theorem~\ref{ACI}, suppose in addition that the residue field of $S$ is infinite and that
$I$ is  generically a complete intersection. Let $J$ be an ideal generated by $g := \codim I$ general elements
of $I$, and let $K := J:I$. For all $\rho$ the module $I^\rho (S/K)$ is an $\omega$-self-dual Cohen-Macaulay
$S/K$-module. In particular, Conjecture~$\ref{main conjecture0}$ is true under these additional hypotheses.
\end{corollary}

\begin{proof}
 We may assume that $I$ is not a complete intersection.
Thus by \cite{Cowsik-Nori} the analytic spread of $I$ is $g+1$, and $K$ is a geometric link of $I$. By Theorem~\ref{ACI}, the ring $S/K$ is Gorenstein, and $I(S/K)$ is generated by
a nonzerodivisor. The conclusion is now immediate.
\end{proof}

\smallskip


When $J:I$ is a geometric $(g+1)$-residual intersection, $I/J$ itself has good properties:

\begin{proposition}
  Let $S$ be a Gorenstein local ring with infinite residue field and let $I\subset S$ be generically a complete intersection of codimension $g$
  such that
 $S/I$ is Cohen-Macaulay.
Let $J\subset I$ be generated by $g+1$ general elements of $I$ and set $K = J:I$.
The module $I/J$ is  $\omega_{S/K}$-self dual and is a Cohen-Macaulay module of dimension $\dim S -g-1 = \dim S/K$.
\end{proposition}

\begin{proof}
 We note that the ideal $K = J:I$ has codimension $\geq g+1$, hence is a $g+1$-residual
intersection of $I$. Since $S/I$ is Cohen-Macaulay and $I$ is generically a complete intersection, $K$ has
codimension exactly $g+1$ (\cite[Proposition 1.7a]{U}). A result of van Straten and Warmt implies that $I/J$ is $\omega_{S/K}$-self dual; see Theorem 2.1 of~\cite{EU} where Huneke's simplified proof is given.

Let $J_{g}\subset J$ be the ideal generated by $g$ general elements. We obviously have $J_{g}:J\supset J_{g}:I \supset K(J_{g}:J)$. 
Every 
associated prime of $J_{g}:I$ has codimension $g$, and hence does not contain $K$. Thus, $J_{g}:J = J_{g}:I $. Therefore,
$J/J_{g} \cong S/(J_{g}:I)$, which has has depth $\dim S -g$. It follows that $\depth S/J\geq \dim S - g-1$, so 
$\depth I/J\geq \dim S -g -1$; that is, $I/J$ is a maximal Cohen-Macaulay $S/K$ -module.
\end{proof}

\begin{remark} There are certainly further phenomena to explain in these directions. For example, let
$I$ be the ideal of $2\times 2$ minors of the generic $3\times 3$ matrix over the polynomial ring $S$ in 
9 variables. We have $\codim I= 4$ and
$\ell(I) = 9$ because $I$ is of
linear type (\cite[Theorem 2.4]{H3}). 

For $s$ with $\codim I = 4 \leq s\leq 8=\ell(I)-1$, let $K_s = J_s:I$ and $R_s = S/K_s$, where $J_s$ is generated by $s$ general forms of degree 2 in $I$. 

The rings $S/I^\rho$ have depth 0  and linear resolution for all $\rho\geq 2$ \cite[Theorem 5.1]{R} but the modules $IR_s$ are maximal Cohen-Macaulay $R_s$-modules
and:
\begin{itemize}
\item $\depth R_4 = 5$; so this ring is a Cohen-Macaulay almost complete intersection;
\item $\depth R_5 = 1,\ \depth R_6 = 0$
\item $R_7$ and $R_8$ are Gorenstein rings of dimensions $2$ and 1, respectively.
\end{itemize}
\end{remark}

\vspace{.03in}
\section{Necessity of the hypotheses}\label{S examples}
Some examples showing that the hypotheses in Conjecture~\ref{main conjecture0} cannot simply be dropped. The following
examples were discovered and checked using the program Macaulay2 \cite{M2}.

\begin{example}
 We first the ideal $K$ of a smooth rational quartic in $\PP^3_\QQ$ as a general link:
 Let $K' \subset \QQ[x_1,\dots,x_4]$ be the ideal of the smooth rational quartic in $\P^3_\QQ$, and
 let $J'\subset K'$ be the ideal generated by two general cubic forms in $K'$. Let
 $I' = J':K'$, which is the ideal of a smooth genus 1 quintic curve in $\P^3$. It turns out that
 $I'$ is minimally generated by 5 cubic forms. If $a$ is a general cubic in $I'$ and  $I := (J',a)$,
 then $I$ is generated by 3 forms of degree 3, and $\ell(I) = 3$.
Finally, let $J$ be the ideal generated by two general cubics in $I$. 

The ideal $K = J:I$ is again the ideal of a smooth rational quartic, and thus neither $R:=S/K$ nor any power of the
principal ideal $IR$ is Cohen-Macaulay. 

Here all the assumptions of Conjecture~\ref{main conjecture0} are satisfied except that  $I$ is not unmixed.
\end{example}

\begin{example}
Let $k$ be an infinite field, and let $X\subset \P^{d-1}_k$ be an abelian surface embedded by a complete linear series of high degree. Let 
$S = k[x_1,\dots,x_d]$ be the homogeneous coordinate ring of $\P^{d-1}$, and let $I_X$ be the homogeneous ideal of $X$. The canonical module $\omega$ of $S/I_X$ is isomorphic to $S/I_X$ as a graded module, and $S/I_X$ is not Cohen-Macaulay because $\HH^1(\cO_X) \neq 0.$ Let $I$ be a homogeneous geometric link of $I_X$, so that $I$ is an unmixed but not Cohen-Macaulay almost complete intersection that is generically a complete intersection.

Let $K$ be any homogeneous link of $I$ with respect to a subset of a system of homogeneous minimal generators, chosen sufficiently generally that $K$ is a geometric link of $I$. Since $I$ is an almost complete intersection, the ideal $IR$ is generated by a single nonzerodivisor. The canonical module of $R:=S/K$ is isomorphic to $IR$, up to shift---that is, $R$ a quasi-Gorenstein ring. 

Since $IR$ is generated by a nonzerodivisor and $R$ is not Cohen-Macaulay, no power of $IR$ can be a Cohen-Macaulay module (though all powers of $IR$ are $\omega_R$-self-dual). 

Here all the assumptions of Conjecture~\ref{main conjecture0} are satisfied except possibly that $I$ is  generated in a single degree.

Now specialize to the case where $X$ is the Segre embedding of the product of two smooth cubic curves in $\PP^2_\QQ$. In fact $I$ is not generated in a single degree.
\end{example}

\begin{example}
Let  $S = \QQ[x_1, \dots, x_7]$ and let
$$
I = (x_1x_4x_7^4,x_5x_6^2x_7^3,x_1x_4x_5^2x_6^2,x_1^2x_3x_4x_5x_6).
$$
The ideal $I$ has codimension 2 and analytic spread 4.
If $J$ is generated by 3 general
forms of degree 6 in $I$, then $K = J:I$ is a 3-residual intersection.
 Because $IR$ is principal,
the high powers of $IR$ are isomorphic, up to a shift, to $S/(J:I^\infty)$, which is not Cohen-Macaulay.
It is interesting to note that $R = S/K$ is Cohen-Macaulay.

Here all the assumptions of Conjecture~\ref{main conjecture0} are satisfied except that I has 
embedded components and  $K$ is not 
a geometric residual intersection of $I$.
 \end{example}

The following examples show that none of the  hypotheses listed in Conjecture~\ref{analytic spread 1} and
Question~\ref{analytic spread 1 with truncation} 
can simply be dropped. 

%

\begin{example}
 Let $H\subset \QQ[x_1,x_2,x_3]$ be the ideal of maximal 
minors of the matrix 
$$\begin{pmatrix}
 x_1^2 & x_1^2 & x_2^2x_3\\
 x_2^2 & x_2^2 & x_1x_2x_3\\
x_2^2 & x_3^2 &x_1x_2 x_3\\
 x_3^2 & 0 & x_1^3
\end{pmatrix}\ .
$$
Let $R$ be the ring defined by the link of $H$ with respect  to the minors deleting the first and second rows.
The ring $R$ is Cohen-Macaulay and generically a complete intersection of dimension 1. The canonical ideal $I$ is 
generated in a single degree, but no power of $I$ is self-dual. 
(Note that because $R$ is 1-dimensional only powers up to the reduction number of $I$ need to be checked.)

Here all the assumptions of Conjecture~\ref{analytic spread 1} are satisfied except that $R$ is not reduced.
\end{example}

\begin{example}
Let  $R$ be the homogeneous coordinate
ring of 11 points in $\PP^2_\QQ$, 6 of which are general and 5 are on a line. 
The ring $R$ is reduced and 1-dimensional, but the canonical ideal $I$ has no self-dual power.
Here all the assumptions of Conjecture~\ref{analytic spread 1} are satisfied except that
 $I$ is not equigenerated. 
 
In this case the fractional ideal $I$ is generated in degrees $-3$ and $-2$; If we take $I'$ to be
the truncation of $I$ in degree -2, then the square of $I'$ is self dual, as is every higher power, giving
a positive answer for Question~\ref{analytic spread 1 with truncation}. 
\end{example}

\begin{example}\label{5points}
Let $R$ be the 
homogeneous coordinate ring of  5 points in $\PP^{2}_\QQ$, of which 3 lie on one line and three on another line (the point of intersection
is one of the 5 points):
$$
[1:0:0], [0:1:0], [0:0:1], [1:1:1], [1:1:0].
$$ 
The ideal  
$$
I = (x_1^2+x_2^2+x_3^2, x_2x_3)R
$$
is equigenerated, but has no self-dual power.

Here all the assumptions of Question~\ref{analytic spread 1 with truncation} are satisfied except that $I$ and the canonical ideal have no power in common.

Curiously, the minimal $R$-free resolutions of $I$ and the canonical ideal of $R$ have the same graded betti numbers for at least 10 steps. However, $I$ and $I^2$ are both generated by 2 elements, whereas the square and cube of the canonical ideal require 3 generators.
\end{example}

\begin{example}\label{no power}
 Let $S = \QQ[x_1,\dots,x_7]$ and let
$$
I = (x_3x_5x_7^4,\ x_2^2x_6^2x_7^2,\ x_2x_3x_4x_5x_6^2,\ x_1x_2x_3x_4x_5x_6).
$$
The codimension of $I$ is 2, its analytic spread is 4, and $I$ satisfies $G_4$. 

If $J$ is generated by 3 general
forms of degree 6 in $I$, then $K = J:I$ is a geometric $3$-residual intersection, necessarily of codimension 3, as $3<4 = \ell(I)$,
and $K$ is reduced by Bertini's Theorem.
The ring $R = S/K$ is Cohen-Macaulay but not Gorenstein. 

Since $IR$ is principal, generated
by a nonzerodivisor, no power of $IR$ can be $\omega_R$-self-dual. 

Here, as in Example~\ref{5points}, all the assumptions of Question~\ref{analytic spread 1 with truncation} are satisfied except that $IR$ and the canonical ideal have no power in common.
\end{example}

\bigskip

\vbox{\noindent Author Addresses:\par
\smallskip
\noindent{David Eisenbud}\par
\noindent{Mathematical Sciences Research Institute,
Berkeley, CA 94720, USA}\par
\noindent{de@msri.org}\par

\medskip
\noindent{Craig Huneke}\par
\noindent{Department of Mathematics, University of Virginia, Charlottesville, VA 22904, USA}\par
\noindent{huneke@virginia.edu}\par

\medskip
\noindent{Bernd Ulrich}\par
\noindent{Department of Mathematics, Purdue University, West Lafayette, IN 47907, USA}\par
\noindent{ulrich@math.purdue.edu}\par

}

\end{document}